\documentclass[10pt,twoside]{article}
\usepackage{amsmath,amscd,amssymb, amsbsy,euscript}
\usepackage{syntonly}
\usepackage{Latex-document}
\input epsf
\newcommand{\comment}[1]{}
\newcommand{\bq}{\begin{equation}}
\newcommand{\tq}{\end{equation}}

\newcommand{\bs}{\begin{split}}
\newcommand{\es}{\end{split}}

\def\al{\alpha}
\def\be{\beta}
\def\de{\delta}

\def\ga{\gamma}
\def\ep{\varepsilon}
\def\ka{\kappa}
\def\th{\theta}

\def\l{\ell}

\def\Om{\Omega}
\def\si{\sigma}

\def\cC{\mathcal C}

\def\cE{\mathcal E}
\def\cG{\mathcal G}
\def\cV{\mathcal V}

\def\limn{\lim_{n \to \infty}}
\def\limk{\lim_{k \to \infty}}
\numberwithin{equation}{section}

\newtheorem{thm}{Theorem}

\newtheorem{Def}{Definition}

\newtheorem{con}{}

\def\lbl(#1){}
\def\EQN(#1){\lbl(#1)
\protect\begin{equation}\protect\label
{#1}}

\def\CON(#1){\lbl(#1) \begin{con}
\label{#1} \end{con}}
\renewcommand{\thesection}{\arabic{section}}
\renewcommand{\thesubsection}{\thesection.\arabic{subsection}.}
\def\Section#1{\section{\hskip -1em . \hskip 0.6em#1}}

\def\volumeno{I}

\title{\bf Some Highlights of Percolation \vskip 6mm}
\author{Harry Kesten\thanks{Department of Mathematics, Cornell University, Malott Hall,
Ithaca NY 14853, USA. E-mail:
kesten@math.cornell.edu}\vspace*{-0.5cm}}
\date{\vskip -8mm}

\begin{document}

\maketitle

\thispagestyle{first} \setcounter{page}{345} \vskip -2mm
\begin{abstract}

\vskip 3mm

We describe the percolation model and some of the principal results and
open problems in percolation theory. We also discuss briefly the
spectacular recent progress by Lawler, Schramm, Smirnov and Werner
towards understanding the phase transition of percolation (on the
triangular lattice).

\vskip 4.5mm


\end{abstract}

\vskip 10mm

\Section{Introduction and description of the percolation model} \setzero

\vskip-5mm \hspace{5mm}

Percolation was introduced by Broadbent and Hammersley (see [14],[15])
as a probabilistic model for
the flow of fluid or a gas through a random medium.
It is one of the simplest models which has a phase transition, and is
therefore a valuable tool for probabilists and statistical physicists
in the study of phase transitions.
For many
mathematicians percolation on general graphs may be of interest
because it exhibits relations between probabilistic and topological
properties of graphs. On the applied side, percolation has been used
to model the spread of a disease or fire, the spread of rumors or
messages, to model the displacement of oil by water,
to estimate whether one can build nondefective integrated
circuits with certain wiring restrictions.


We shall give a brief survey of some
of the important results
obtained for this model
\comment
{, including some of the recent results by
Lawler, Schramm, Smirnov and Werner on conformal invariance and
powerlaws. Several}
and list some open problems.
The present article
is only a very restricted survey and its references (in particular to
the physics literature)
are far from complete. We apologize to the authors of relevant articles which
we have not cited.
Earlier surveys are in [42], [21], [22], [37], [38],
and the reader can find more
elaborate treatments in the books [20],
[63]
and [55].

The oldest (indirect) reference to percolation that I know of is a
problem submitted to the Amer. Math. Monthly (vol 1, 1894, pp. 211-212)
in 1894 by De Volson Wood, Professor of Mechanical Engineering at the
Stevens Inst. of Technology in Hoboken NJ. Here is the text of the problem.

\medskip
``An actual case suggested the following:

\smallskip
An equal number of white and black balls of equal size are thrown into
a rectangular box, what is the probability that there will be
contiguous contact of white balls from one end of the box to the
opposite end ? As a special example, suppose there are 30 balls in the
length of the box, 10 in the width and 5 (or 10) layers deep.''

\medskip
Even though percolation theory was not invented to answer this
problem, it naturally came to study problems of this kind. By the
way, we still have no answer to De Volson Wood's problem. Percolation
as a mathematical theory was invented by Broadbent and Hammersley
([14],[15]). Broadbent wanted to model the spread of a gas or
fluid through a random medium of small channels which might or might
not let gas or fluid pass. To model these channels he took the edges
between nearest neighbors on $\mathbb Z^d$ and made all edges
independently {\it open} (or passable) with probability $p$ or {\it
closed} (or blocked) with probability $1-p$. Write $P_p$ for the
corresponding probability measure on the configurations of open and
closed edges (with the obvious $\si$-algebra generated by the sets
determined by the states of finitely many edges).
A {\it path} on $\mathbb Z^d$ will be a sequence (finite or infinite)
$v_1, v_2, \dots$ of vertices of $\mathbb Z^d$ such that for all $i \ge 1,\;
v_i$ and $v_{i+1}$ are adjacent on $\mathbb Z^d$.
The edges of such a path are the edges $\{v_i,v_{i+1}\}$
between successive vertices
and a path is called $\it open$ if all its
edges are open. Broadbent's original question amounted to asking for
\begin{equation}\label{1.1}
P_p\{\exists \text{ an open path on $\mathbb Z^d$ form $\bold 0$ to
}\infty \}.
\end{equation}
This question has an obvious analogue on any infinite connected graph
$\cG$
with edge set $\cE$ and
vertex set $\cV$.
Again one makes all edges independently open or closed with
probability $p$ and $1-p$, respectively, and one denotes the corresponding
measure on the edge configurations by $P_p$. $E_p$ is
expectation with respect to $P_p$. An open path is defined
as before with $\cG$ taking the role of $\mathbb Z^d$. A path
$(v_1, v_2,\dots)$ is called {\it self-avoiding}
if $v_i \ne v_j$ for $i \ne j$.
\eqref{1.1} now is replaced by
\begin{equation}\label{1.2}
P_p\{\exists \text{ an infinite self-avoiding open path
starting at $v$}\},
\end{equation}
with $v$ any vertex in $\cV$.

The preceding model is called {\it bond-percolation}.
There is also an analogous model, called {\it site-percolation}. In
the latter model all edges are assumed passable, but the vertices are
independently open or closed with probability $p$ or $1-p$,
respectively. An open path is now a path all of whose vertices are
open. One is still interested in \eqref{1.2}. Site percolation is more
general than bond percolation in the sense that the positivity of
\eqref{1.2} for some $v$ in
bond-percolation on a graph $\cG$ is equivalent to the positivity of
\eqref{1.2} for some $v$ in site-percolation on
the covering graph or line graph of $\cG$. However, site percolation on a
graph may not be equivalent to bond percolation on another graph (see
[40], Section 2.5 and Proposition 3.1).

{\it Unless otherwise stated we restrict ourselves in the remaining sections
to site percolation}.
We shall often use $\cV$ and $\cE$ to
denote the vertex and edge set of whatever graph we are discussing at
that moment,
without formally introducing the graph as $\cG=(\cV, \cE)$.
It should be clear from the context what $\cV$ and $\cE$ stand for in
such cases. For $A \subset
\cV$, we shall use
$|A|$ to denote the number of vertices in $A$.
Further if $A,B$ and $C$ are sets
of vertices, then $A \leftrightarrow B$ means that there exists an
open path from some vertex in $A$ to some vertex in $B$, while
$A \overset C  \leftrightarrow B$ means that there exists an
open path with all its vertices in $C$,
from some vertex in $A$ to some vertex in $B$. In particular,
with some abuse of notation, we have
\[
\{|\cC(v)| = \infty\}=\{v \leftrightarrow \infty\}.
\]

\begin{Def} We call a graph $\cG = (\cV, \cE)$ {\it quasi-transitive}
if there is a finite set of vertices  $V_0$, such that for each vertex
$v$ there is a graph automorphism of $\cG$ which maps $v$ to one of
the vertices in $V_0$.
\end{Def}
All vertices which can be mapped by a graph automorphism to a
fixed $v_0 \in V_0$ are equivalent for our purposes. In  a
quasi-transitive graph each vertex is equivalent to one of finitely
many vertices. A special subclass is formed by the {\it transitive graphs},
which have $|V_0| = 1$, so that all vertices are equivalent
for our purposes. (For example,
the Cayley graph of a finitely generated group is transitive.)

We shall restrict ourselves here to
graphs which are
\begin{equation}\label{1.3a}
\text{connected, infinite but locally finite,
and quasi-transitive}.
\end{equation}
\comment{
The graph $\cG$ is (vertex) {\it transitive}
if for any $v',v'' \in \cV$ there is a graph
automorphism which takes $v'$ to $v''$.  Thus all vertices of
such a graph are equivalent for our purposes.
\endcomment}
Graphs which satisfy \eqref{1.3a}
automatically have countable vertex sets and edgesets.
We define, for $v \in \cV$,
\begin{equation}\label{1.3}
\bs
\th^v(p) &= P_p\{v \leftrightarrow \infty\}\\
&=P_p\{\exists \text{ an infinite self-avoiding open path
starting at $v$}\}.
\end{split}
\end{equation}
For a quasi-transitive graph
$\th^v(p) = \th^{v_0}(p)$ for some $v_0 \in V_0$.
It is an easy consequence of the
FKG inequality that either $\th^v(p) > 0$ for all $v$ or $\th^v(p) =
0$ for all $v$ (see [40], Section 4.1).
We call $\th^v(p)$ the {\it percolation
probability (from $v$)}. Much of the earlier work on percolation theory deals
with properties of the function $p \mapsto \th^v(p)$, or more generally
with the full distribution of the so-called cluster sizes. The {\it cluster}
$\cC(v)$ of the vertex $v$ is the set of all points which are
connected to the origin by an open path. By convention, this always
contains the vertex $v$ itself (even if $v$ itself is closed
in the case of site percolation).
\comment{
In the case of bond percolation
the clusters are the maximal
components of the graph generated by the open edges, or
single vertices all of whose incident edges are closed. Similarly, in
the case of site percolation on $\cG$,
}
The clusters are the maximal components
of the graph with vertex set $\cV$ and with an edge between two sites
only if they are adjacent on $\cG$ and are both open. $\th^v(p)$ is just
the $P_p$-probability that $|\cC(v)|= \infty$.

\Section{Existence of phase transition and related properties of the
critical probability} \setzero

\vskip-5mm \hspace{5mm}

The most important property of the percolation model is that it exhibits
a {\it phase transition}, that is, there exists a threshold value $p_c$
such that the global behavior of the system is quite different in the
two regions $p < p_c$ and $p > p_c$. To make this more precise let us
consider the percolation probability as a function of $p$. It is
non-decreasing. This is easiest seen from Hammersley's ([31]) joint
construction of percolation systems for all $p \in [0,1]$ on $\cG$. Let
$\{U(v), v \in \cV\}$ be independent uniform $[0,1]$ random variables.
Declare $v$ to be $p$-open if $U(v) \le p$. Then the configuration of
$p$-open vertices has distribution $P_p$ for each $p \in [0,1]$. Clearly
the collection of $p$-open vertices is nondecreasing in $p$ and hence
also $\th(\cdot)$ is nondecreasing. Clearly $\th^v(0) = 0$ and
$\th^v(1)= 1$. Roughly speaking the graph of $\th^v(\cdot)$ (for a fixed
$v$) therefore looks as in figure 1, but not all the features exhibited
in this figure have been proven.

\epsfverbosetrue \epsfxsize=200pt \epsfysize=200pt \centerline
{\epsfbox{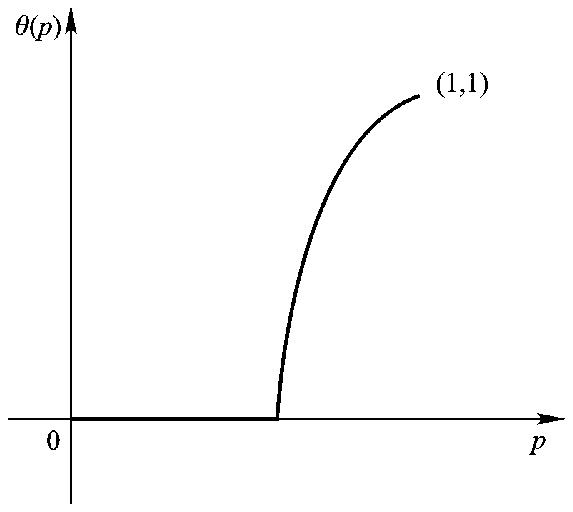}}

\centerline {Figure 1: Graph of $\th$. Many aspects of this graph
are still conjectural.}

\bigskip
\noindent
The {\it critical probability} is defined as
\begin{equation}\label{1.4}
p_c = p_c(\cG) = \sup\{p:\th^v(p)=0\}.
\end{equation}
As remarked after \eqref{1.3} this is independent of $v$.
By definition we then have
\[
P_p\{|\cC(v)|= \infty\} = 0 \text{ for } p< p_c, v \in \cV,
\]
so that
\begin{equation}\label{1.6}
\text{all clusters are finite a.s. $[P_p]$ when $p < p_c$}.
\end{equation}
On the
other hand, for $p > p_c$ there is a strictly positive $P_p$-probability
that $|\cC(v)|$ is infinite. It then follows from Kolmogorov's
zero-one law that
\begin{equation}\label{1.6a}
P_p\{\text{some }|\cC(v)| = \infty\} = 1,\quad p > p_c.
\end{equation}
Thus the global behavior of the system is quite different for $0 \le p <
p_c$ and for $p_c < p <1$. We therefore can say that there is a phase
transition at $p_c$, provided the intervals $[0,p_c)$ and $(p_c,1]$ are
both nonempty. It is easy to see from a so-called Peierls argument (just
as in [29]) that $p_c(\cG) > 0$ for any graph $\cG$ of bounded degree
(and hence certainly if \eqref{1.3a} holds). It is much harder to do
show that $p_c(\cG)< 1$ holds for certain $\cG$. Hammersley [30] proved
this for bond-percolation on $\mathbb Z^d$, but a similar argument works
for site-percolation and various other periodic graphs. (Basically, we
say that $\cG$ is {\it periodic} or can be {\it periodically imbedded
in} $\mathbb R^d$ if $\cV$ can be imbedded in $\mathbb R^d$ (with $d \ge
2$) such that $\cV$ as well as the edges of $\cG$, as represented by the
straight line segments between the pairs of vertices adjacent in $\cG$,
form a subset of $\mathbb R^d$ which is invariant under translations by
$d$ linearly independent vectors. If this is the case we call $d$ the
dimension of $\cG$. We refer the reader to [40], Section 2.1 for
details.) Thus\vskip -2mm
\begin{thm}
\[
0 < p_c(\mathbb Z^d) < 1.
\]
\end{thm}

Thus, at least on $\mathbb Z^d$, there really is a phase transition.
On any graph one says that the system is in the
{\it subcritical (supercritical) phase} if $p < p_c$ (respectively, $p
> p_c$).
Because percolation is such a simple model with a phase transition,
percolation has received a great deal of attention from physicists.
Percolation is one of the Potts models,
corresponding to the parameter $q$ in the Potts model equal to 1; the
famous Ising model for magnetism is essentially the same as
the Potts model with
$q=2$. One hopes that understanding of the percolation model
will help understand all the Potts models and even the more general
Fortuin-Kasteleyn or random cluster models (see [23], which also
explains the relation, due to Fortuin and Kasteleyn,
between random cluster models and Potts models).

The exact value of $p_c(\cG)$ is known only for a handful of graphs,
and all of these are periodic two-dimensional graphs.
This leads to
\newline
{\bf Open problem 1:} Find $p_c(\cG)$ for a wide class of graphs.
\newline
However, it is generally agreed that the solution to this problem
would not have any explanatory value. The critical probabilities which
have been determined so far depend heavily on special symmetry
properties of the underlying graph, and the values of these critical
probabilities vary with the graph. One has therefore moved on to
properties which are believed to be shared by large classes of graphs;
see Section 4 below. The rigorously
known critical probabilities can be found in [38], Chapter 3.
Here we merely mention the one case which will be important
later on:
\begin{equation}\label{2.10}
p_c(\text{site percolation on triangular lattice}) = \frac 12.
\end{equation}
Also known is the following asymptotic result, both for the site and
for the bond version:
\begin{equation}\label{2.11}
p_c(\mathbb Z^d) \sim \frac 1{2d} \text{ as }d \to \infty.
\end{equation}
This has been proven by several people; [35] gives the best higher
order terms in \eqref{2.11}.

One can define another critical probability as the threshhold value
for the finiteness of the clustersize of a fixed vertex.
Thus,
\begin{equation}\label{2.12}
p_T(\cG) = \sup\big \{p: E_p\{|\cC(v)|\} = \infty \big\}.
\end{equation}
Since $P_p\{|\cC(v)| = \infty \} > 0$ for $p > p_c$, it is obvious that
$E_p\{|\cC(v)|\} = \infty$ for all $p > p_c$, so that $p_T(\cG) \le
p_c(\cG)$. It was a crucial step in establishing the known values for
$p_c$ to show that $p_T(\cG) = p_c(\cG)$. The original proof of this
fact was only for bond percolation on $\mathbb Z^2$ ([39]; this proof
made strong use of crossing probabilities similar to those appearing in
De Volson Wood's problem in Section 1). Proofs of $p_T = p_c$ for some
other special lattices are in [65] and [40]). Later Menshikov ([51]) and
Aizenman and Barsky ([1]) gave independent and different proofs of
exponential decay of the distribution of $|\cC(v)|$ for $p < p_c$. This
is a cornerstone of the subject and is of course a much stronger
statement than $p_T = p_c$.\vskip -2mm
\begin{thm}(Menshikov and Aizenman and Barsky) Assume that
$\cG$ is periodic. Then for $p < p_c(\cG)$
 there exists constants $0 < C_1,C_2 < \infty$ such that
\begin{equation}\label{2.12a}
P_p\{|\cC(v)| \ge n\} \le C_1 e^{-C_2n},\quad n \ge 0.
\end{equation}
\end{thm}

\eqref{2.12a} gives a basic estimate  for the subcritical phase.
By an earlier ``subadditivity'' argument of [45] \eqref{2.12a} can be
sharpened to a ``local limit theorem'' (see [20], Theorem 6.78):
for each $p < p_c$
there exists a $0 < C_3(p) < \infty$ such that
\begin{equation}\label{2.12b}
\limn -\frac 1n \log P_p\{|C(v)| = n\} = C_3(p).
\end{equation}

These results give us a measure of control over the subcritical phase.
In the supercritical phase many estimates rely on another fundamental
result of percolation theory, which was proven by Grimmett and Marstrand
[24]. The simplest form of the result is as follows:\vskip -2mm
\begin{thm}\label{GM}
\begin{equation}\label{2.13} p_c(\mathbb Z^d) = \limk p_c(\mathbb Z_+^2 \times \{1, 2,
\dots,k\}^{d-2}).
\end{equation}
One may replace $\mathbb Z_+^2$ by $\mathbb Z^2$ here.
\end{thm}

The graph appearing in the right hand side here consists of a finite
number of copies of the first quadrant in $\mathbb Z^2$  or of the
whole $\mathbb Z^2$. Thus (before the limit is taken)
this graph looks very much like $\mathbb Z^2$ and many of the special
tools for  percolation on $\mathbb Z^2$ can be applied to this graph.
Because of this one could prove a number of results on $\mathbb Z^d$
for $p > \limk p_c(\mathbb Z_+^2 \times \{1, 2,
\dots,k\}^{d-2})$. Theorem \ref{GM} now shows that these results hold
throughout the supercritical regime (at least when $\cG = \mathbb Z^d$
or a similar graph). As an example of this situation we mention a
result of [43], namely the right hand inequality in \eqref{2.14}
(the left hand inequality is due to [3]):
For site percolation on $\mathbb Z^d$ with $p > p_c(\mathbb Z^d)$
there exist $0< C_4(p), C_5(p) < \infty$ such that
\begin{equation}\label{2.14}
C_4(p) \le - \frac 1{n^{(d-1)/d}} \log P_p\{|C(v)| =  n\} \le C_5
\end{equation}
for all large $n$.
\newline
{\bf Open problem 2:} Does
\[
\limn - n^{-(d-1)/d}\log P_p\{|C(v)| = n\} \text{ exist ?}
\]
(under the conditions for \eqref{2.14}). The reader should notice
the contrast between \eqref{2.12}, \eqref{2.12a} --- which give
exponential decay for the clustersize distribution in the
subcritical case --- and \eqref{2.14} which corresponds to a
``stretched exponential'' for the tail of the clustersize in the
supercritical case. The tail of this distribution at criticality,
i.e., for $p = p_c$ will be discussed in Section 4.

\Section{Uniqueness of infinite clusters and properties of the
percolation probability} \setzero

\vskip-5mm \hspace{5mm}

It is  natural to ask ``how many infinite
clusters can there be ?'' In  [52] it is shown that for
periodic graphs for each $p$, exactly one of the
following three situations prevails:
\newline
$P_p\{\text{there is no infinite open cluster}\}=1$,
\newline
$P_p\{\text{there is exactly one infinite open cluster}\}=1$ or
\newline
$P_p\{\text{there are infinitely many infinite open clusters}\}=1$.
\newline
As pointed out in [58], the proof of [52] carries
over to  any
quasi-transitive graph by a zero-one law for events which are
invariant under graph automorphisms.
Of course, the first alternative here holds for $p < p_c$, but can the last
situation occur for some $p \ge p_c$ ? The first proof that this is
impossible on $\mathbb Z^d$ is in [4]. This proof was improved and
generalized a few times, but the most elegant, and by now standard,
proof is due to Burton and Keane [16]. Their method works for any
amenable graph. To make this precise we define for any set
$W \subset \cV$,
\[
\partial W = \{w \in \cV: w \notin W \text{ but $w$ is adjacent to some
$v \in W$}\}.
\]
We call the graph $\cG$ {\it amenable} if there exists a sequence
$\{W_n\} \subset \cV$ for which $|\partial W_n|/|W_n| \to 0$.\vskip -2mm
\begin{thm}(Burton and Keane) If $\cG$ satisfies \eqref{1.3a},
and if $\cG$ is amenable, then for all $p \in [0,1]$
\begin{equation}\label{2.15}
P_p\{\text{there exist more than one infinite open cluster}\} = 0.
\end{equation}
\end{thm}
The proof of this result is the same as in [16], except that one
should argue on the {\it expected number of encounter points} where
Burton and Keane use the ergodic theorem to make the number of
encounter points itself large. (We owe this observation to
O. H\"aggstr\"om.)
Simple examples (such as a regular tree) show that \eqref{2.15} does
not have to hold for nonamenable graphs. This is one example of
a relation between percolation properties and algebraic/topological
properties of the underlying graph (see [10], [50] and [8] and
some of their
references for other examples).
What can be said about uniqueness/nonuniqueness  in the nonamenable
case ? Benjamini and Schramm [10] introduced a further critical
probability:
\begin{equation}\label{2.18}
p_u = p_u(\cG) :=\inf \{p:\text{a.s. $[P_p]$ there is a unique infinite
cluster}\}.
\end{equation}

By definition $p_u \ge  p_c$. We have $p_c < p_u =1$ on a regular
$b$-ary tree (in which all vertices have degree $b+1$) with $b \ge 2$.
The first example of a graph with $p_c < p_u < 1$ was given in [25].
Note that there is no a priori reason why uniqueness should be monotone
in $p$, that is why uniqueness a.s. $[P_{p'}]$ should imply uniqueness
a.s. $[P_{p''}]$ whenever $p'' \ge p'$. This has been proven to be the
case for graphs satisfying \eqref{1.3a}. More precisely, the following
theorem (and somewhat more) is proven in [57] (see also [26] and
[27]):\vskip -2mm
\begin{thm} Let $\cG$ satisfy \eqref{1.3a} and let
the percolation configurations on $\cG$ be constructed
simultaneously for all $p \in [0,1]$ by Hammersley's  method described
 in the beginning of Section 2. Let $N(p)$ be the number of $p$-open
 infinite clusters. Then a.s.,
\[
N(p) = \begin{cases} 0 &\text{ for } p \in [0,p_c)\\
\infty &\text{ for } p \in (p_c,p_u)\\
1 &\text{ for }p \in (p_u,1].
\end{cases}
\]
\end{thm}
Note that this theorem does not give the value of $N(p)$ at $p=p_c$ or
$p_u$ (see also the lines after Theorem 1.2 in [26] and Open problem
3 below).

Other obvious questions concern the smoothness of the function
$\th^v(\cdot)$, and in particular whether this function is
continuous. Clearly $p \mapsto \th^v(p)$ is always continuous for $p <
p_c$, since $\th^v(p) = 0$ for all such $p$. Russo [53]
noted that $\th^v(\cdot)$ is everywhere right continuous and [11]
proved that (under \eqref{1.3a}) if for some fixed
$p_0 > p_c$ there is a.s. $[P_{p_0}]$ a unique infinite cluster, then
$\th^v(\cdot)$ is also left continuous at $p_0$. Thus, under
\eqref{1.3a}, the remaining
problem is
\newline
{\bf Open problem 3:} Is $p \mapsto \th^v(p)$ (left) continuous on
$[p_c,p_u]$ ?
\newline
On $\mathbb Z^d$ which has $p_c(\mathbb Z^d) = p_u(\mathbb Z^d)$,
continuity is equivalent to $\th(p_c)= 0$. It has long been
conjectured that this is the case. It is known that this holds for
$d \ge 19$ by the theory of Hara and Slade [34]; actually this
deals with bond percolation, but should go through also for site
percolation on $\mathbb Z^d$. It also follows from work of Harris
[36] and the author [40], Theorem 3.1, that continuity holds when
$d=2$ (both for bond and site percolation). [8] and [9] prove that
on a Cayley graph of a non-amenable group there is no percolation
at $p_c$.

\Section{Behavior at and near $p_c$} \setzero

\vskip-5mm \hspace{5mm}

{\it From now on we shall restrict ourselves to transitive graphs which
are periodically imbedded in $\Bbb R^d$, so that the origin is a vertex
of the graph}. Since all vertices are equivalent in a transitive graph,
we drop the superscript $v$ from various quantities such as $\th(p)$; we
further write $\cC$ for the open cluster of the origin.

We saw in \eqref{2.12a}  and \eqref{2.14} that the probability of
a cluster of size $n < \infty$ decays exponentially or as a stretched
exponential in the subcriticial and supercritical regime,
respectively. The behavior at criticality is quite different. In fact,
it is believed that there exists constants
$0 < C_i< \infty$ such that
\begin{equation}\label{4.1}
C_6n^{-(d-1)/2} \le P_{p_c}\{|C(v)| \ge n \} \le C_7n^{-C_8}.
\end{equation}
Indeed, for periodic graphs in dimension $d=2$ the left hand
inequality
is proven in [12], but the argument remains valid in any dimension.
In an Abelian sense one knows even more (see the proof of Proposition
10.29 in [20]).
The right hand inequality of \eqref{4.1}
for certain two-dimensional graphs can be found in
[40], Theorem 8.2, while for $d \ge 19$ with $\de = 2$ it
follows from [6] and [33].
It is natural to conjecture that
\begin{equation}\label{4.2}
P_{p_c}\{|C(v)| \ge n \} \approx n^{-1/\de}
\end{equation}
for some $\de = \de(\cG)> 0$, where $a(n) \approx b(n)$ means
$\log(a(n)/\log b(n) \to 1$ as $n \to \infty$. (One may conjecture
that \eqref{4.2} and similar relations below hold with an  even
stronger interpretation of $\approx$ but we shall not pursue this
here.)

\eqref{4.2} is one example of a so-called {\it power law}. Another
(conjectured) power law is for $P_{p_c} \{ v' \leftrightarrow v''\}$.
It is believed that for some constant $\eta$
\begin{equation}\label{4.3}
P_{p_c} \{ v' \leftrightarrow v''\} \approx |v'-v''|^{2-d-\eta}
\text{ as } |v'-v''| \to \infty.
\end{equation}
Here $|v|$ denotes the $\l^1$ norm of the image of $v$ under the
imbedding into $\mathbb R^d$.
Again this is supported by the following partial result for periodic
graphs whose image under the periodic imbedding into $\Bbb R^d$
is invariant under permutations of the coordinates.
For such graphs \eqref{4.1} implies that
there exist constants $0 < C_i < \infty$ such that
\begin{equation}\label{4.3a}
C_9|v'-v''|^{-C_{10}} \le P_{p_c} \{ v' \leftrightarrow v''\} \le
C_{11}|v'-v''|^{-C_{12}}.
\end{equation}

Physicists also conjectured that various quantities behave like powers
of $|p-p_c|$ as $p \to p_c, p \ne p_c$.
Such conjectures are analogues of results which were known (often on a
nonrigorous basis) or conjectured for related models. They were also
rigorously known for quite some time for percolation on regular trees.
In addition, Hara and Slade in a series of important papers (see in
particular [33], [34]) have developed the so-called lace expansion
technique to give us a good understanding of percolation in high
dimensions. Roughly speaking, they prove many of the physicists
conjectures for bond percolation
on $\mathbb Z^d$ with $d \ge 19$, by showing that
most quantities show mean field
behavior near $p_c$, that is, they have the same singularity on
$\mathbb Z^d$ with $d \ge 19$ as on a regular tree. It is believed
that this will remain true for $d > 6$. In fact Hara and Slade can
prove their results for percolation in any dimension $> 6$ for what they
call ``spread out'' models. (These have $\mathbb Z^d$ as vertex set
but there may be some open bonds between points which are not nearest
neighbors on $\mathbb Z^d$.)

The most common of the
conjectured power laws (with the traditional names for the exponents)
are as follows. Here $A(p) \approx B(p)$ means $\log A(p)/\log B(p)
\to 1$ as $p \to p_c$ (with $p \ne p_c$).
\begin{equation}\label{4.5}
\Big (\frac d{dp}\Big)^3 \Big[\sum_{n=1}^\infty \frac 1n P_p\{|\cC| = n\}\Big]
= \Big (\frac d{dp}\Big)^3 E_p \{|\cC|^{-1}\} \approx |p-p_c|^{-1-\alpha},
\end{equation}
\begin{equation}\label{4.6}
\th(p) \approx (p-p_c)^\be, \quad p \downarrow p_c,
\end{equation}
\begin{equation}\label{4.7}
\chi(p) := E\{|\cC(v)|; |\cC(v)| < \infty\} \approx |p-p_c|^{-\ga}.
\end{equation}
Another power law is supposed to hold for
the so-called {\it correlation length}, $\xi(p)$. Intuitively
speaking, if  $p \ne p_c$, the correlation length is the minimal size a
cube should have so that one can detect from a typical percolation
configuration in such a cube that $p$ is not equal to $p_c$.
On scales which are small with respect to the correlation
length, the system is expected to behave as if it is critical. On the
other hand,
on scales which are large with respect to the correlation length, one
should be able to partition the system into cubes of edgelength equal
to a large multiple of the correlation length and regard these
cubes as ``supersites'';
for $p >p_c$ (respectively $p <p_c$) these supersites should behave as
sites in site percolation with a $p$ value close to 1 (respectively
close to 0). On such scales the details of the lattice other than its
dimension should play little role, if any. Several
possible formal definitions are in use for the correlation
length. Here we define the correlation  length $\xi(p)$ by
\begin{equation}\label{4.8}
[\xi(p)]^{-1} = \lim_{n \to \infty} -\frac 1{n} \log
P_p\{\bold 0 \leftrightarrow ne_1, \bold 0 \nleftrightarrow
\infty\},
\end{equation}
where $e_1$ is the first coordinate vector. Strictly speaking, [19]
only proves that this is a good definition for (bond or site) percolation
on $\mathbb Z^d$, but this definition should make sense with minor
changes for percolation on general periodic graphs.
The conjectured powerlaw then takes the form
\begin{equation}\label{4.9}
\xi(p) \approx |p-p_c|^{-\nu}.
\end{equation}
Other power laws have been conjectured for electrical conductance and
for the graph-theoretical length of an open crossing between opposite
faces of a cube.

In all these cases proofs of power {\it bounds} instead of actual power laws
are known for many graphs which are periodically imbedded in $\mathbb R^d$
with $d=2$ or $d$ large (see [40], Chapter 8, [20], Chapter 10,
[33], [34]). For instance,
\begin{equation}\label{4.10}
C_{13}|p-p_c|^{-1} \le \chi(p) \le C_{14}|p-p_c|^{-C_{15}} \text { for
} p < p_c.
\end{equation}
In fact, the left hand inequality holds for all $d$ (see [5], [20],
Theorem 10.28).

Most remarkable is the conjecture of ``universality''. That is, it is
generally believed that each of the so-called critical
exponents $\al, \be, \ga, \de, \eta, \nu$
depends for periodic graphs on the dimension $d$ only, and not
on the details of the graph $\cG$. For
instance, they should have the same value for bond and for site
percolation on $\mathbb Z^2$ and on the triangular lattice. This is in
contrast to the critical probablity $p_c$, which definitely
{\it does} depend on the details of $\cG$. For this reason the
principal concern these days is to establish power laws and
universality, and little attention is being paid to open problem
1. (See next section for more on what is now known.)

There are also nonrigorous arguments to derive simple relations
between various of these exponents. These are the so-called scaling
laws:
\begin{equation}\label{4.10a}
\al + \be (\de+1) = 2
\end{equation}
\begin{equation}\label{4.11}
\ga + 2\be = \be (\de + 1)
\end{equation}
\begin{equation}\label{4.12}
\ga = \nu(2-\eta)
\end{equation}
\begin{equation}\label{4.13}
d\nu = \ga + 2\be \text{ for } 2 \le d \le 6.
\end{equation}
The last relation, which involves the dimension $d$
is called a hyper-scaling law. These scaling relations, except \eqref{4.10a},
 have been established for  many graphs
with $d=2$ (assuming that the exponents exist; see [41])
and are also known
for high $d$. There is also a ``conditional proof'' of the
hyperscaling relation . That is,
\eqref{4.13} (or rather the relation $(2-\eta) = d(\de -1)/(\de+1)$)
has been shown to be implied by other
(not yet established) laws for percolation ([13]). It is widely
believed that these other laws hold for $3 \le d \le 6$.

There are predictions by physicists of the values of these exponents
when $d = 2$ or when $d$ is large. In fact, as we already pointed out,
 it is believed that all
these exponents are even independent of $d$ for $d > 6$. The
existence of these exponents and their
predicted values have now been proven to be correct when $\cG$ is the
triangular lattice ([60], [61], [62], [48]). It is also known that these
exponents (except for $\al$ and perhaps $\eta$)
exist and take the same values as on a regular tree
for $d \ge 19$ ([33], [34], [32]).

 This section raises the obvious and very extensive
\newline
{\bf Open problem 4:} Prove power laws, universality and scaling relations.

  In the next section we shall describe some of the progress made on
this problem in dimension 2. We already mentioned the work of Hara and
Slade in high dimensions. No progress has been made
in dimensions 3, 4 and 5. So we may pose a more modest problem for
these dimensions.
\newline
{\bf Open problem 5:} Find upper and lower power bounds for $\th(p), \chi(p)$,
$\xi(p)$ and $P_{p_c}\{|\cC| \ge  n\}$ when $3\le   d \le  5$.

As we pointed out above, bounds on one side are already known for
most of these quantities, but as far as we know no bound of the form
\begin{equation}\label{4.15}
             \xi(p) \le  C_{16}|p - p_c|^{-C_{17}}
\end{equation}
has been proven for $3 \le  d \le  5$, not even for $p < p_c$.
This is probably
the most fundamental bound to prove, from which several other bounds
might follow. Note that it is not hard to see that on $\Bbb Z^d$
\begin{equation}\label{4.16}
\xi(p) \ge  C_{18}|\log( p_c - p)|^{-(d-1)/d}(p_c - p)^{-1/d}
\text{ for }p < p_c.
\end{equation}
Indeed, the proof of the left hand inequality in \eqref{4.3a} actually gives
\begin{equation}\label{4.17}
P_{p_c}\{0 \overset {[0,n]^d} \leftrightarrow ne_1\} \ge  C_9n^{-3(d-1)}.
\end{equation}
From this one trivially has for $p < p_c$
\begin{equation}\label{4.18}
\bs
P_p\{\bold 0 \leftrightarrow  kne_1\} &\ge \big[P_p\{\bold 0
\overset {[0,n]^d} \leftrightarrow
ne_1\} \big]^k\\
&\ge \big[\big(\frac p{p_c}\big)^{(n+1)^d}
P_{p_c}\{\bold 0 \overset {[0,n]^d}\leftrightarrow ne_1\}\big]^k\\
&\ge \big[\big(\frac p{p_c}\big)^{(n+1)^d}C_9n^{-3(d-1)}\big]^k.
\end{split}
\end{equation}
Now take $n = \big[|\log(p_c - p)|\big]^{1/d}(p_c - p)^{-1/d}$
and estimate $\xi(p)$ from
\[
\bs
[\xi(p)]^{-1}&= \limk -\frac 1{kn} \log P_p\{ \bold 0 \leftrightarrow kne_1,
\bold 0 \nleftrightarrow \infty\}\\
&= \limk - \frac 1{kn}  \log P_p\{ \bold 0 \leftrightarrow kne_1\}
\text{ for } p <
p_c.
\end{split}
\]

A somewhat different aspect of the behavior of critical percolation
concerns the random variable
\begin{equation}\label{4.19}
N(v) := \inf\{\text{number of closed vertices in any path from $\bold
0$ to $v$}\}.
\end{equation}
If no percolation occurs for $p = p_c$, then $N(v) \to \infty$ as
$v \to \infty$, a.s. $[P_{p_c}]$. For bond percolation on
$\mathbb Z^2$ it is known that $[\si_{p_c}(N(v))]^{-1}[N(v)- E_{p_c}N(v)]$
satisfies a central limit theorem with $E_{p_c}N(v) \asymp
\log|v|$ and $\si_{p_c}(N(v))$, the standard deviation of $N(v)$,
of order $[\log |v|]^{1/2}$ (see [44]). On
$\mathbb Z^d$ with $d \ge 3$ it is only known ([18]) that $N(v) =
O(|v|^\ep)$ a.s. $[P_{p_c}]$, for every $\ep > 0$.
\newline
{\bf Open problem 6:} Improve the bound for $N(v)$ and find a limit theorem
for $N(v)$ in dimension $\ge  3$.

\Section{Conformal invariance and SLE} \setzero

\vskip-5mm \hspace{5mm}

{\it In this section we only consider graphs which are periodically
imbedded in $\mathbb R^2$}.  Special attention will be paid to site
percolation on the triangular lattice.

We already briefly discussed the
interpretation of the correlation length in the preceding
section. In view of \eqref{4.3a}, the definition \eqref{4.8} assigns
the value $\infty$
to the correlation length when $p=p_c$,
at least if $P_{p_c}\{|\cC| = \infty\} =
0$, as is widely believed (and is known for $d = 2$ or $d \ge 19$).
Thus the correlation length is not a useful length scale for critical
percolation. Other than the spacing between vertices, there seems to
be no lengthscale which plays a role for critical percolation.
In this case one may hope to take some sort of limit without normalization
of a critical percolation system in a larger and larger region. It is not
clear in what topology one should take a limit. Matters look somewhat
friendlier if one fixes a region and considers a limit as the spacing
between vertices tends to zero. Even then it is not clear what topology
will be most useful for taking a limit. A discussion of these issues can be
found in [2] and the beginning of [59]. Putting this problem aside,
let us first ask for limits of simple quantities
such as crossing probabilities.  Let $D$ be a Jordan domain in
$\mathbb R^2$ with
a smooth boundary and let $A_1$ and $A_2$ be two disjoint arcs of
$\partial D$. Identify $\cG$ with its periodic imbedding in $\mathbb
R^2$.
(This imbedding is not unique, but for the present purposes we can
just fix some imbedding.)
We can then define $\de \cG$ as the result of multiplying the image of
$\cG$ under
the imbedding by a factor $\de > 0$.  This image has vertices located at
$\{ \de v : v \in \cV\}$ and edges between two points
$\de v', \de v''$ if and only if $v'$
and $v''$ are adjacent in $\cG$.  For any percolation configuration on
$\cG$ we say that there exists an open path on $\de \cG$ from $A_1$ to
$A_2$ in $D$ if there
is an open path $v_1,\dots,v_m$  on $\cG$ such that $\de v_i \in D$
for $2 \le  i \le  m - 1$
and the edge between $\de v_1$ and $\de v_2$ intersects $A_1$ and
the edge between $\de v_{m-1}$  and $\de v_m$  intersects $A_2$.
We then define
\begin{equation}\label{5.1}
h(D, A_1,A_2, \de):=P_{p_c}\{\exists \text{ open path on $\de \cG$ from
$A_1$ to $A_2$ in $D$}\},
\end{equation}
and ask whether this has a limit as $\de \downarrow 0$.
(Here is were contact is made
with De Volson Wood's problem in the Amer. Math. Monthly.) It is
conjectured that this limit, call it $h(D, A_1, A_2)$, exists, and moreover
that it is conformally invariant. By this we mean that if $\phi$
is a conformal map from $D$ onto $D' = \phi(D)$ which extends to a
homeomorphism between $\overline D  := $ closure of $D$ and $\overline
D'$, then
\begin{equation}\label{5.2}
h(D, A_1, A_2) = h(\phi(D), \phi(A_1), \phi(A_2)).
\end{equation}
Conformal invariance of a limit of critical percolation had been
conjectured by physicists (see [17] and its references) on the grounds that
this had been found in related models.  The stress on studying this
for crossing probabilities is due to [46], which also credits Aizenman
with the formulation of conformal invariance for crossing
probabilities
(actually in a slightly more general form than \eqref{5.2}).  Cardy used
conformal invariance and the Riemann mapping theorem to equate
$h(D, A_1, A_2)$ to $h(\mathbb H, [z, 0], [1, \infty))$,
where $\mathbb H$ is the upper half plane
and $z \in (-\infty, 0)$ a suitable point on the boundary of $\mathbb
H$, i.e., the real axis.  He then derived (nonrigorously) a
differential
equation for $h(\mathbb H, [z, 0], [1, \infty))$ as a function of
$z$. From this he obtained an explicit
formula for $h(\mathbb H, [z, 0], [1, \infty))$, and hence for
$h(D, A_1, A_2)$ in special cases, such as when D is a rectangle and
$A_1, A_2$ two opposite sides of $D$.
In an astonishing paper Smirnov [60] succeeded in showing that for
site percolation on the triangular lattice, the limit $h(D, A_1, A_2)$
indeed exists and is conformally invariant. To do this Smirnov introduces an
extra variable $z \in \overline D$ , and considers
\[
\bs
f(z, D, A_1, A_2, \de) &:= P_{p_c}\{\exists \text{ self-avoiding
open path on $\de\cG$
from }A_1\\
&\text{to $B_1 \cup  A_2$ in $D$ which separates $z$ from }B_2\},
\end{split}
\]
where $B_1, B_2$ are the arcs on $\partial D$ between $A_1$ and $A_2$
(i.e., the boundary of $D$ consists of the four arcs $A_1, B_1, A_2,
B_2$
and one successively traverses these arcs as one goes around the
boundary of $D$ in one direction).  He now shows that any limit of
$f(z, D, A_1, A_2, \de)$ along a
subsequence $\de_n \downarrow 0$ is a harmonic function of
$z\in D$ which has to satisfy certain boundary conditions which
uniquely determine the limit.
Therefore $\lim_{\de \downarrow 0} f(z, D, A_1, A_2, \de_n)$ exists.
Moreover, the limit is conformally invariant, because it is
characterized
as the harmonic function which satisfies a certain boundary
condition.
The original problem for the crossing probabilities $h(D, A_1, A_2,
\de)$ can be treated as a special case, by
letting $z$ approach the single point in $\overline A_2 \cap \overline
B_1$. One can find the limit function $h(D, A_1,A_2)$ explicitly if
$D$ is a rectangle, and $A_1, B_1, A_2, B_2$ its sides, and
thereby one can recover Cardy's formula.

Somewhat before Smirnov, Schramm [59] had introduced {\it stochastic
Loewner evolutions} (SLE) in order to describe a scaling limit of
growing random sets (and in particular the scaling limit of loop erased
random walk in dimension two, and related processes). For percolation,
the simplest version of SLE is probably the so-called chordal SLE
(see [54]), described as follows. Let $\mathbb H$ and
$\overline{\mathbb H}$  be the open upper and
closed upper half plane, respectively, and let $\{B(t)\}_{t\ge  0}$
be a standard Brownian motion starting at 0.
Let $g_t(z)$ be the solution of the Loewner equation
\begin{equation}\label{5.3}
\frac {\partial g_t(z)}{\partial t} = \frac 2{g_t(z)-\xi(t)},
\quad g_0(z) = z
\end{equation}
with $\xi(t) = \sqrt \ka B(t)$ for some parameter $\ka > 0$.
The solution to \eqref{5.3} exists for $t < \tau (z) :=
\inf\big\{s : 0 \text{ is a limit point of the }$ set $
\{g_u(z)-\xi(u), u < s\}\big\}$. Define
$H_t := \{z \in \mathbb H :\tau(z) > t\}, K_t = \{z \in \overline H
:\tau (z) \le  t\}$. Chordal SLE$_\ka$ is the collection of maps $\{g_t
: t \ge  0\}$. It turns out that $g_t$ is the unique conformal
homeomorphism from $H_t$ onto $\mathbb H$ for which
$\lim_{z \to \infty} [g_t(z) - z] = 0$. It is shown in [54] that for
$\ka \ne 8$ there exists a continuous path $\ga : [0, \infty) \to
\overline {\mathbb H}$  such that $K_t$ is the hull of $\ga[0, t]$,
that is, $K_t$ is the closure of the union of the bounded components of
$\overline {\mathbb H}\setminus  \ga[0, t]$. In many situations one
can also start with the path $\ga$ and
then define $g_t$ as the conformal homeomorphism from its hull $K_t$
onto $\mathbb H$. This must then satisfy a Loewner equation \eqref{5.3}.
$\ga$ is called the {\it trace} of the corresponding
SLE process.

Schramm ([59]) showed that the scaling limit of loop erased random
walk can be described by an analogue of SLE$_2$ in the unit disc.
([59] still had to assume that this scaling limit exists and is
conformally invariant, but this has since been proven in [49]).  In [59]
Schramm expresses the belief that SLE$_6$ is
appropriate for the description of the scaling limit of the boundary
of percolation clusters. This has been proven to be correct for
percolation on the triangular lattice.

\bigskip
\epsfverbosetrue \epsfxsize=330pt \epsfysize=230pt \centerline
{\epsfbox{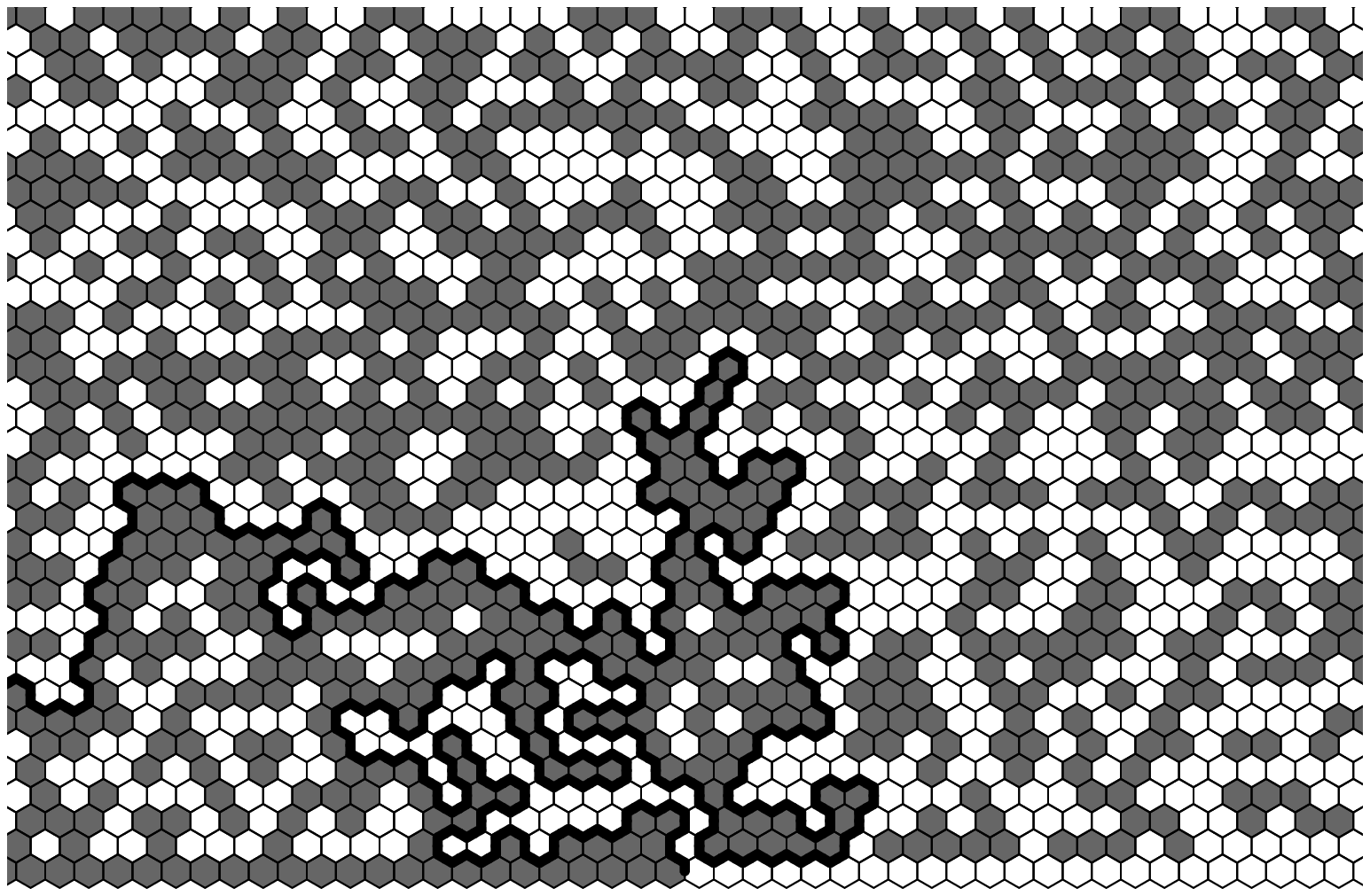}}

\begin{center}
\begin{minipage}[h]{10cm}
Figure 2: The exploration process, which separates the open (white)
hexagons from the closed (black) ones. We thank Oded Schramm for
providing us with this figure.
\end{minipage}
\end{center}

\bigskip
\noindent
To give a specific example, consider
the hexagonal lattice, imbedded in $\mathbb R^2$ in such a way that
the hexagonal faces which intersect the x-axis have their centers on
this axis and that the origin lies on
the common boundary of two such faces.
Make the hexagonal faces in the
upper half plane independently open or closed with probability
1/2. If one thinks of the centers of the hexagonal faces as vertices of
the triangular lattice then one sees that this is equivalent to critical
site percolation on the triangular
lattice on $\mathbb H$
(recall that its critical probability equals 1/2).
Now impose the boundary
 condition that all faces
with center on the positive (negative) x-axis are open
(closed, respectively). There is then a curve,$\ga_\de$,
in the upper half plane and running on
the boundaries of some of the hexagonal faces, which starts at 0 and
traverses the boundary between the open cluster of the positive x-axis and
the closed cluster of
the negative x-axis. This curve is called
the {\it exploration process}; see Figure 2.
The distribution of this curve $\ga_\de$ converges to the distribution
of the trace of SLE$_6$ (as the mesh size goes to zero,
and using the Hausdorff metric on the space of curves, determined up
to parametrization)
(see [60], [61]). Actually these references discuss the analogous
situation on an equilateral triangle instead of $\mathbb H$ and
concentrate on showing the existence of the limit. The
identification of the limit  as SLE$_6$ is based on the work of
Lawler, Schramm and Werner ([47], [64]).

\comment
{This result was proven by
showing that any subsequential limit of
$\ga$ is conformally invariant and has a certain ``restriction
property" which states
that in some sense $\ga$ does not feel the boundary. (Note that one can
construct $\ga$ from ``local'' information only; at any step $\ga$
turns to the right (left) if its tip has a closed (respectively, open)
hexagon in front of it.) Moreover, some probabilities similar to
crossing probabilities are the same for the limit curve
as for SLE$_6$, by virtue of
Cardy's formula. These properties uniquely determine the distribution
of the limit of the exploration process and show that is the same as SLE$_6$.
}

To prove this result Smirnov ([60], [61])
first uses a compactness argument to show that any sequence $\de_n
\downarrow 0$ has a subsequence along which the distribution of
$\ga_{\de_n}$ converges to some distribution $\mu^{e.p.}$ on H\"older
continuous curves. Then he proves that $\mu^{e.p.}$ is independent of
the subsequence $\{\de_n\}$ by showing that $\mu^{e.p.}$ has certain
properties which characterize SLE$_6$. This of course also shows that
$\mu^{e.p.}$ is the distribution of the trace of
SLE$_6$. The second step relies on a reduction of various
$\mu^{e.p.}$-probabilities to crossing probabilities of the form
\eqref{5.1} and on the existence and conformal invariance of the
limit of \eqref{5.1}. In addition it relies on a ``locality
property.'' Note that one can
construct $\ga_\de$ from ``local'' information only; at any step $\ga_\de$
turns to the right (left) if its tip has a closed (respectively, open)
hexagon in front of it.

\comment
{
In addition one needs (an analogue
of) the exploration process on a simply connected domain
$\Om$, from a point $a$ to a point $b$ on the boundary of $\Om$, and
the conformal invariance of its limit distribution.

 This implies that the analogue
of the weak limit --- develops ``without feeling the boundary.''
Moreover its law is a conformal invariantof $(\Om, a,b)$.

a strong relationship between the
distribution of the exploration process on $\Bbb H$ and its analogue
on a conformal image of $\Bbb H$.
}

SLE turns out to be the perfect tool for calculating critical exponents.
Lawler, Schramm, Smirnow and Werner in [48] and [62] were able to use
the correspondence with SLE$_6$ to prove for percolation on the
triangular lattice, not only Cardy's formula, but also the power laws
\eqref{4.2}, \eqref{4.3}, \eqref{4.6}, \eqref{4.7} and \eqref{4.9} with
the values for  $\nu, \be, \ga, \de$ and $\eta$ which were predicted by
physicists (see [62] for relevant references). It was further shown in
[7] by Beffara that the Hausdorff dimension of the trace of SLE$_6$ is
7/4. Thus, this is also the Hausdorff dimension of the exploration
process ``in the scaling limit.'' This dimension had already been
predicted by Saleur and Duplantier [56].

\label{lastpage}

\end{document}